\documentclass{article}
\usepackage[utf8]{inputenc}
\usepackage[english]{babel}

\usepackage{indentfirst} % to indent the first paragraph after a section tag

\usepackage[ruled,vlined]{algorithm2e} % creat algorithms in Latex

\usepackage{amsthm}
\theoremstyle{remark} % The command sets the styling right below it.
\newtheorem*{remark}{Remark}
\theoremstyle{definition} % The command sets the styling right below it.
\newtheorem{example}{Example}

\usepackage{amsmath}
\theoremstyle{plain} % change the style of the output of the environments
\newtheorem{theorem}{Theorem}
\newtheorem{lemma}{Lemma}
\newtheorem{corollary}{Corollary}

\usepackage{amsfonts,stmaryrd,amssymb} % Math packages

\usepackage[
backend=biber,
style=ieee,
dashed=false,
sorting=none,
maxnames=9
]{biblatex} % a modern program to process bibliography information
\usepackage{csquotes}
\DeclareFieldFormat[article, inbook, incollection, inproceedings, misc, thesis, unpublished]{title}{#1} % 'fix' the quotes around the title of the paper
\renewbibmacro{in:}{} % to removes the 'in'
\newtoggle{bbx:datemissing} % year at the end
\renewbibmacro*{date}{\toggletrue{bbx:datemissing}}
\renewbibmacro*{issue+date}{
  \toggletrue{bbx:datemissing}
  \iffieldundef{issue}
    {}
    {\printtext[parens]{
       \printfield{issue}}}
  \newunit}
\newbibmacro*{date:print}{
  \togglefalse{bbx:datemissing}
  \printdate}
\renewbibmacro*{chapter+pages}{
  \printfield{chapter}
  \setunit{\bibpagespunct}
  \printfield{pages}
  \newunit
  \usebibmacro{date:print}
  \newunit}
\renewbibmacro*{note+pages}{
  \printfield{note}
  \setunit{\bibpagespunct}
  \printfield{pages}
  \newunit
  \usebibmacro{date:print}
  \newunit}
\renewbibmacro*{addendum+pubstate}{
  \iftoggle{bbx:datemissing}
    {\usebibmacro{date:print}}
    {}
  \printfield{addendum}
  \newunit\newblock
  \printfield{pubstate}}

\usepackage{dsfont} % display the characteristic function

\usepackage{enumerate}

\usepackage[framemethod=tikz]{mdframed} % Allows defining custom boxed/framed environments

\usepackage{geometry} % Required for adjusting page dimensions and margins

\usepackage[utf8]{inputenc} % Required for inputting international characters
\usepackage[T1]{fontenc} % Output font encoding for international characters

\usepackage{XCharter} % Use the XCharter fonts
\usepackage[title]{appendix}

\mdfdefinestyle{commandline}{
	leftmargin=10pt,
	rightmargin=10pt,
	innerleftmargin=15pt,
	middlelinecolor=black!50!white,
	middlelinewidth=2pt,
	frametitlerule=false,
	backgroundcolor=black!5!white,
	frametitle={Command Line},
	frametitlefont={\normalfont\sffamily\color{white}\hspace{-1em}},
	frametitlebackgroundcolor=black!50!white,
	nobreak,
}

\title{Convergence of Gradient Algorithms for Nonconvex $C^{1+\alpha}$ Cost Functions\footnote{This work was partially supported by National Key R$\&$D Program of China (2018YFA0703900) and National Natural Science Foundation of China (Grant Numbers 11631004 and 12031009).}} % Title of the assignment

\author{Zixuan Wang\footnote{Department of Finance and Control Sciences, Shanghai Center for Mathematical Sciences, Fudan University, Shanghai, 200433, China (email: 17110840010@fudan.edu.cn)}\,, \qquad Shanjian Tang\footnote{Department of Finance and Control Sciences, School of Mathematical Sciences, Fudan University, Shanghai, 200433, China (email: sjtang@fudan.edu.cn)}} % Author name and email address

\date{} % University, school and/or department name(s)
\begin{document}

\allowdisplaybreaks[1]
\setlength{\abovedisplayskip}{4pt}   % reset space above displayed equations
\setlength{\belowdisplayskip}{4pt}   % reset space below displayed equations

\maketitle % Print the title

\setlength{\parskip}{0.2em}
%----------------------------------------------------------------------------------------
%	ABSTRACT
%----------------------------------------------------------------------------------------
\noindent\rule{\textwidth}{1pt}
\vspace{-2.25em}
\section*{Abstract}
\vspace{-0.25em}
\noindent This paper is concerned with convergence of stochastic gradient algorithms with momentum terms in the nonconvex setting.
A class of stochastic momentum methods, including stochastic gradient descent, heavy ball, and Nesterov’s accelerated gradient, is analyzed in a general framework under mild assumptions.
Based on the convergence result of expected gradients, we prove the almost sure convergence by a detailed discussion of the effects of momentum and the number of upcrossings.
It is worth noting that there are not additional restrictions imposed on the objective function and stepsize.
Another improvement over previous results is that the existing Lipschitz condition of the gradient is relaxed into the condition of H{\"o}lder continuity.
As a byproduct, we apply a localization procedure to extend our results to stochastic stepsizes.\\
%\hspace*{\fill}
\noindent\rule{\textwidth}{1pt}

%----------------------------------------------------------------------------------------
%	INTRODUCTION
%----------------------------------------------------------------------------------------
\section{Introduction}
In recent years, deep learning has witnessed an impressive string of empirical success in the area of image classification, speech recognition, natural language processing, etc.
Due to the rapid growth in the scale of modern datasets, finding the minimum of a function $f$ with an iterative procedure has become very popular, especially to minimize the training error of deep networks.
We study the classical unconstrained stochastic programming problem of the form
\begin{equation}
\min_{\mathbf{x}} \ f(\mathbf{x}) \ \hat{=} \ \mathbb{E}^\mathbb{P}[h(\mathbf{x}, Z)],
\end{equation}
where $\mathbf{x} \in \mathbb{R}^d$, $h$ is a measurable function, and $Z$ is a random element with a known or unknown probability law $\mathbb{P}$.
When a discrete distribution is considered, i.e., $Z$ represents a random index obeying the uniform distribution on the finite set $\{z_1, \dots, z_n\}$, the stochastic optimization problem is given in the form
\begin{equation}
\min_{\mathbf{x}} \ f(\mathbf{x}) \ \hat{=} \ \frac{1}{n}\sum_{i=1}^n h(\mathbf{x}, z_i).
\end{equation}
For instance, in supervised learning, each $z_i$ represents a training sample, $\mathbf{x}$ is the parameter of the model, and $f$ represents the training loss.
The landscape of the cost function remains blurred in the case $f$ does not have special structure.
Hence, sometimes we relax the problem to finding a critical point of $f$.
\par
A heuristic approach is to take steps proportional to the negative of the gradient, which is known as \emph{Gradient Descent} (GD).
It was first set forth in \cite{cauchy1847} by Cauchy dated 1847.
The convergence was proved in \cite{curry1944} for the least squares problem where each component function $h(\cdot, z)$ is the square of a continuously differentiable function.
Then related research progressed towards convex programming problems and new variants with momentum were proposed.
In a celebrated work \cite{polyak1964}, Polyak came up with the \emph{Heavy Ball} (HB) method, i.e., GD with exponentially weighted memory, to speed up the convergence for convex optimization.
Furthermore, in \cite{nesterov1983}, \emph{Nesterov's Accelerated Gradient} (NAG) method achieved the optimal convergence rate for a convex function $f$ with Lipschitz continuous gradient.
The application of quasi-Newton methods is also explored.
For example, Becker and LeCun \cite{beckerlecun1989} used diagonal approximations of the Hessian matrix.
See also \cite{polyak1987}, \cite{polyakjuditsky1992}, \cite{nesterov2004}, and \cite{nesterov2012} for classical results regarding GD.
\par
The exact value of the gradient is required for all the aforementioned algorithms.
Nevertheless, the effective computation of the gradient is too cost in a large-scale optimization problem.
Sometimes our access to $f$ or $\nabla f$ is limited when considering simulation-based problems or problems with unknown $\mathbb{P}$.
Therefore, these deterministic algorithms are restrictive.
To address issues, stochastic gradient algorithms originated from \cite{robbinsmonro1951} and \cite{kieferwolfowitz1952} where the randomized gradient substituted for the exact value.
Not only stochastic methods keep the complexity per iteration constant with respect to the scale of the problem, but they are also likely to escape local minima.
For this reason, stochastic versions of gradient algorithms such as \emph{Stochastic GD} (SGD), \emph{Stochastic HB} (SHB), and \emph{Stochastic NAG} (SNAG) recently have regained interest.
\par
However, there exists a gap between practical success and theoretical explorations.
Nonconvex deep neural networks are usually trained with decaying learning rates while many convergence results are gained for programming problems with convexity (see \cite{nesterov2004}, \cite{polyak1987}, \cite{ghadimietal2015}, \cite{bottouetal2018}, and references therein) or fixed stepsizes (see \cite{lietal2019}, \cite{yanetal2018}, \cite{lessardetal2016}, and references therein).
Although there are some results for nonconvex problems under mild conditions, these works are often performed in a case-by-case manner.
Hence, convergence properties and the complexity are still open in theory, especially for momentum and adaptive methods in a systematic approach.
\par
Here we focus on the sufficient conditions of the convergence in a unified treatment.
Using the method in \cite{ghadimilan2013} and the algorithm framework of \cite{yanetal2018} flexibly, we show $L^2$ convergence of these gradient algorithms in Theorem 1.
This allows us to bound the effects of momentum in the original observation Lemma 3 lying at the center of our demonstration.
Then we prove almost sure convergence in Theorem 2 by the construction of a supermartingale and a detailed discussing of upcrossings, which are extensions of analysis in \cite{bertsekastsitsiklis1996}.
Specifically, the main contributions of this paper are summarized in the following.
\vspace{-4pt}
\begin{itemize}
\item Firstly, we demonstrate almost sure convergence of stochastic momentum methods including SGD, SHB, and SNAG under mild assumptions.
To the best of our knowledge, the theoretical assurance of almost sure convergence of SNAG for a nonconvex $f$ has not been proved.
The majority of previous works in nonconvex setting analyze asymptotic behavior in the sense of distribution or $L^2(\mathbb{P})$.
Nonetheless, in practical application we always complete the training of neural networks few times, which is equivalent to taking several sample points from the corresponding distribution.
Compared with distribution properties, the analysis of a fixed sample path of stochastic algorithms seems more instructive.
\vspace{-2pt}
\item Secondly, the Lipschitz condition of $\nabla f$ is relaxed into the condition of H{\"o}lder-continuity.
Indeed, the derivative of such an uncomplicated $\mathbb{R}^1$ function $f(x)=x^{4/3}\mathds{1}_{|x|\le1} +(2x^{2/3}-1)\mathds{1}_{|x|>1}$ satisfies the $\frac{1}{3}$-H{\"o}lder condition but not Lipschitz condition.
Additionally, we do not need constant or decreasing stepsizes, or coercive condition of $f$, i.e., $\lim_{\Vert\mathbf{x}\Vert \to \infty}f(\mathbf{x}) = +\infty$.
These restrictions are often assumed to hold.
For example, the coercive condition is requisite in the proof of almost sure convergence of SHB in \cite{gadatetal2018}.
\vspace{-2pt}
\item Lastly, stochastic stepsizes are also analyzed and the convergence of some perturbed adaptive algorithms with momentum terms is obtained.
\end{itemize}
\par
\vspace{-2pt}
The rest of the paper is organized as follows.
In Section 2, we survey some recent related works.
After introducing notations as well as the analytical framework, Section 3 gives the main results and proofs.
Section 4 studies stochastic stepsizes and adaptive algorithms.
Finally, Section 5 concludes with a brief discussion.

%----------------------------------------------------------------------------------------
%	main content
%----------------------------------------------------------------------------------------
\section{Related Works}
There are abundant works on the convergence of stochastic gradient algorithms.
We concentrate on articles discussing nonconvex situations.
\par
The first complete proof of almost sure convergence in the absence of convexity is fulfilled by Bertsekas and Tsitsiklis in \cite{bertsekastsitsiklis1996} and \cite{bertsekastsitsiklis2000}, where the classical SGD is considered.
The core is that the process $f(\mathbf{x}_t)$ can be shown to be approximately a supermartingale.
Then the Robbins-Siegmund lemma \cite{robbinssiegmund1971} and a thorough inspection of upcrossings guarantee the almost sure convergence of $f(\mathbf{x}_t)$ and $\nabla f(\mathbf{x}_t)$ respectively.
Since every step is influenced by the lasting memory, their analysis cannot be applied directly to algorithms with momentum terms.
\par
Ghadimi and Lan \cite{ghadimilan2013} demonstrate the convergence of the gradient in $L^2(\mathbb{P})$ sense.
They creatively use a random number $R$ to terminate vanilla SGD.
Then by direct computation of $\mathbb{E} \Vert\nabla f(\mathbf{x}_R)\Vert^2$, the complexity is established.
Along this route, a class of randomized accelerated gradient algorithms including SNAG is proved to converge in $L^2(\mathbb{P})$ sense in \cite{ghadimilan2016}, in which stepsizes must satisfy additional requirements.
Yan et al. \cite{yanetal2018} also propose an enlightening framework to unify SGD, SNAG, and SHB.
However, using the same technique, only expectation convergence analysis is carried out in the case of small constant stepsizes and the bounded gradient.
Their discussion on the influence of momentum is insufficient to obtain further results.
\par
Note that almost sure convergence cannot be derived from these papers.
Notwithstanding the fact that some technical conditions can be relaxed, there is intrinsic difficulty.
Exactly speaking, their demonstration hinges heavily on the termination time $R$.
Given pre-fixed required accuracy $\epsilon$, we can choose suitable $R_\epsilon$ to ensure $\mathbb{E} \Vert \nabla f(\mathbf{x}_{R_\epsilon}) \Vert^2$ is small enough.
When $\epsilon$ tends to zero, the algorithm must terminate at different $R_\epsilon$ which depends on $\epsilon$ and has different distribution on the set of positive integers.
It means that, the approach is unadaptable to the study of asymptotic behavior, especially almost sure convergence.
\par
The ODE approach is introduced by Ljung in \cite{ljung1977}, and extensively developed by \cite{kushnerclark1978}, \cite{kushnershwartz1984}, \cite{fortpages1996}, etc.
The basic idea is approximating discrete-time stochastic algorithms by continuous-time approach where the limit is an ordinary differential equation.
On some stability assumptions, the bridge is built between the behavior of each sample path of the SGD and the related ODE.
With the help of conclusions from the ODE method, almost sure convergence of SHB is obtained in \cite{gadatetal2018} by constructing a Lyapunov function.
The limitation is that the coercive condition of $f$ is required and stepsizes have the form $1/n^r (r \in (0,1])$.
\cite{gitmanetal2019} adopts a different approach based on careful calculation and Levy's extension of the Borel-Cantelli lemma, but only the limit inferior is demonstrated with momentum parameter $\beta_t$ tending to $1$ or $0$.
\par
In spite of limitations, some research is motivated by the ODE approach.
\cite{hoffmanblei2015} uses a multivariate Ornstein-Uhlenbeck process to approximate SGD in the vicinity of a local minimum, but the statement there is in heuristic territory.
The idea is further developed and mathematical aspects are solidified in \cite{lietal2017} and \cite{lietal2019}.
They go beyond OU process approximations and use a class of stochastic differential equations to study the dynamics of SGD, SHB, and SNAG methods.
Nevertheless, these convergence results are established in the weak sense, i.e., convergence in distribution.
\par
Some efforts have been devoted to utilizing control theoretic tools, such as Integral Quadratic Constraints \cite{lessardetal2016}, PID Controllers \cite{anetal2018}, and Regularity Condition \cite{xiongetal2020}, in the analysis of stochastic algorithms with momentum, but they are also under strong assumptions.

\section{Stochastic Gradient Methods with Momentum}
\subsection{Notations and Setup}
In the following, we will write vectors with bold letters.
Let $\Vert \cdot \Vert$ denote the Euclidean norm of a vector and $\langle\cdot, \cdot \rangle$ denote the inner product.
\par
Often, computing the full gradient can be quite expensive.
Furthermore, we cannot gain the exact gradient if the accurate distribution is not known.
Therefore, instead of observing a full gradient of $f$ at $\mathbf{x}$, we assume that we have a \emph{first order oracle} which, given $\mathbf{x} \in \mathbb{R}^d$, returns a noise gradient $\mathbf{g}(\mathbf{x}, \xi)$ where $\mathbf{g}$ is a Borel measurable $\mathbb{R}^d$-valued function and $\xi$ is a random variable on a probability space $(\Omega, \mathcal{F}, \mathbb{P})$.
In the $t$-th iteration of the gradient method, we observe $\mathbf{g}(\mathbf{x}_t, \xi_t)$ and denote it by $\mathbf{g}_t$ for the sake of brevity.
\par
Below, we repeat the \emph{Stochastic Unified Momentum} (SUM) method that deals with nonconvex stochastic programming problems and allows kinds of momentum terms.
This algorithm is firstly proposed in \cite{yanetal2018}.
We replace the constant stepsize in SUM with a variable one.
The next subsection enjoys the convenience provided by this unified framework.\par
\vspace{-2pt}
{\centering
\begin{minipage}{.92\linewidth}
\begin{algorithm}[H]
  \caption{Stochastic Unified Momentum}
  \label{alg1}
  \textbf{Input:} momentum factor $\beta\in[0, 1)$, parameter $s\ge0$, stepsizes $(\gamma_t)_{t\ge0}$, and initial point $\mathbf{x}_0 \in \mathbb{R}^d$.\\
  $\mathtt{1}$\quad Set $t=0$ and $\mathbf{y}_0^s = \mathbf{x}_0$.\\
  $\mathtt{2}$\quad Sample $\xi_t$\ and get $\mathbf{g}_t$.\\
  $\mathtt{3}$\quad Set
  \begin{equation*}
    \mathbf{y}_{t+1} = \mathbf{x}_t - \gamma_t \mathbf{g}_t,
  \end{equation*}
  \begin{equation*}
    \mathbf{y}_{t+1}^s = \mathbf{x}_t - s \gamma_t \mathbf{g}_t.
  \end{equation*}\\
  $\mathtt{4}$\quad Set
  \begin{equation*}
    \mathbf{x}_{t+1} = \mathbf{y}_{t+1} + \beta (\mathbf{y}_{t+1}^s - \mathbf{y}_{t}^s).
  \end{equation*}\\
  $\mathtt{5}$\quad Substitute $t$ with $t+1$ and go to step 2.
\end{algorithm}
\end{minipage}
\par
}
\vspace{2pt}
Setting $s \!=\! \frac{1}{1\!-\!\beta}$, one can easily verify that SUM reduces to SGD with the stepsize $(\frac{\gamma_t}{1\!-\!\beta})$:
\begin{equation}
\mathbf{x}_{t+1} = \mathbf{x}_t - \frac{\gamma_t}{1\!-\!\beta}\mathbf{g}_t.
\end{equation}
When $s=0$, we have SHB:
\begin{equation}
\left\{
\begin{aligned}
\mathbf{m}_{t+1} &= \beta\mathbf{m}_t - \gamma_t\mathbf{g}_t,\\
\mathbf{x}_{t+1} &= \mathbf{x}_t + \mathbf{m}_{t+1},
\end{aligned}
\right.
\end{equation}
whereas when $s=1$, we have SNAG:
\begin{equation}\left\{
\begin{aligned}
\mathbf{y}_{t+1} &= \mathbf{x}_t - \gamma_t\mathbf{g}_t,\\
\mathbf{x}_{t+1} &= \mathbf{y}_{t+1} + \beta(\mathbf{y}_{t+1}-\mathbf{y}_t).
\end{aligned}
\right.
\end{equation}
\par
Throughout this section, we impose the following assumptions on the cost function $f$.
\par \textbf{(A.1)} $f$ is a continuous differentiable function such that
\begin{equation*}
f_* \ \hat{=} \ \inf_{\mathbf{x}}f(\mathbf{x}) > -\infty .
\end{equation*}
\par \textbf{(A.2)} There exists $\alpha \in (0,1]$ such that $f \in C^{1+\alpha}$, i.e., for some $A > 0$,
\begin{equation*}
\Vert\nabla f(\mathbf{x}) - \nabla f(\mathbf{y})\Vert \le A {\Vert\mathbf{x} - \mathbf{y}\Vert}^\alpha, \quad \forall \mathbf{x},\mathbf{y} \in \mathbb{R}^d .
\end{equation*}
\par $f$ is called a \emph{L-smooth} function in the case $\alpha=1$ and $A=L$.
Smoothness assures us that the gradient does not change dramatically within the region around where it is taken, and thus the value of the gradient is informative when gradient algorithms are applied with small stepsizes.
\par We will also need some assumptions on the stochastic oracle $\mathbf{g}(\mathbf{x}, \xi)$.
The expected direction is assumed to be parallel with the gradient, and a bound on the mean square is essential.
\par \textbf{(A.3)} The oracle is an unbiased estimator for the gradient, i.e.,
\begin{equation*}
\mathbb{E}[\mathbf{g}(\mathbf{x},\xi)] = \nabla f(\mathbf{x}), \quad \forall \mathbf{x} \in \mathbb{R}^d .
\end{equation*}
\par \textbf{(A.4)} There exist positive constants $\sigma^2$ and $C$, such that
\begin{equation*}
\mathbb{E}[{\Vert \mathbf{g}(\mathbf{x},\xi) - \nabla f(\mathbf{x}) \Vert}^2] \le \sigma^2 + C {\Vert \nabla f(\mathbf{x}) \Vert}^2, \quad \forall \mathbf{x} \in \mathbb{R}^d .
\end{equation*}
\par \textbf{(A.5)} The random variables $\{\xi_t\}_{t \ge 0}$ are independent of each other.
\par The following stepsize $(\gamma_t)_{t \ge 0}$ is considered.
\par \textbf{(A.6)} $(\gamma_t)_{t \ge 0}$ is a deterministic and nonnegative sequence such that
\begin{equation*}
\sum_{t=0}^\infty \gamma_t = \infty, \quad
\sum_{t=0}^\infty {(\gamma_t)}^{1+\alpha} < \infty .
\end{equation*}
\par
Define $\mathcal{F}_0 \ \hat{=} \ \left\{ \Omega, \varnothing \right\}$ and $\mathcal{F}_t \ \hat{=} \ \sigma(\xi_0, \xi_1, \dots, \xi_{t-1})$ for $t \ge 1$.
Since stepsizes are assumed to be deterministic, we have $\mathbf{x}_t \in \mathcal{F}_t$.
Notice that (A.5) implies the independence of $\mathbf{x}_t$ and $\xi_t$.
Therefore, $\mathbb{E}[{\Vert \mathbf{g}_t - \nabla f(\mathbf{x}_t) \Vert}^2 |\mathcal{F}_t] \le \sigma^2 + C {\Vert \nabla f(\mathbf{x_t}) \Vert}^2$ a.s.
\par
The above assumptions are standard and reasonable for general stochastic algorithms.
For instance, the cost function $f$ of a multi-layer network with the sigmoid activation function $S(x) = \frac{1}{1+exp(-x)}$ and the quadratic loss is L-smooth.
The existing theoretical results are almost based on assumptions (A.1)-(A.6).
See \cite{bertsekastsitsiklis1996}, \cite{bertsekastsitsiklis2000}, and \cite{ghadimilan2013} for SGD, \cite{nesterov2004} and \cite{ghadimilan2016} for SNAG, \cite{gadatetal2018} for SHB, and \cite{lietal2017}, \cite{lietal2019}, \cite{yanetal2018}, and \cite{gitmanetal2019} for unified treatment.
In fact, the coercive condition of $f$ and the boundedness of $\nabla f$ and $\mathbf{g}$ are always assumed in previous works.
\par
To show that the assumption (A.2) is more general, we give several examples of machine learning, in which the cost function satisfies (A.2) but not L-smooth condition.
\begin{example}
We consider the linear \emph{support-vector machine} (SVM) for binary classification.
We are given a training dataset of $n$ points: $(\mathbf{x}_1, y_1), \dots, (\mathbf{x}_n, y_n)$, where each $\mathbf{x}_i$ is a point in $\mathbb{R}^d$ and $y_i \in \{1, -1\}$ indicates the class to which $\mathbf{x}_i$ belongs.
Our aim is to minimize the empirical risk for the smooth hinge loss introduced by \cite{rennie2005}:
\par
\begin{equation}
h_\alpha(v) = \begin{cases}
\frac{\alpha}{\alpha+1} -v, & v \le 0;\\
\frac{1}{\alpha+1}v^{\alpha+1} -v +\frac{\alpha}{\alpha+1}, & 0 < v < 1;\\
0, & v \ge 1,
\end{cases}
\end{equation}
where $\alpha > 0$.
So the cost function is
\begin{equation}
f_\alpha(\mathbf{w}, b) = \frac{1}{n} \sum_{i=1}^n h_\alpha(y_i (\langle\mathbf{w},\mathbf{x}_i\rangle +b)) ,
\end{equation}
where $(\mathbf{w}, b) \in \mathbb{R}^d\times\mathbb{R}$.
When $\alpha \ge 1$, the gradient of $f_\alpha$ is Lipschitz continuous, whereas when $0< \alpha < 1$, $f_\alpha \in C^{1+\alpha}$ but it is not L-smooth.
Note that $h_\alpha$ converges uniformly to the original hinge loss $h(v) = \max \{0, 1-v\}$ as $\alpha \to \infty$. Additionally, the gradient of the corresponding cost function for the hinge loss $h$ is not continuous.
\end{example}
A nonconvex $C^{1+\alpha}$ example can be found in \cite{sypherdetal2019}, where the sigmoid activation function combined with ``$\alpha$-loss'' is considered.
Some $C^{1+\alpha}$ regularization terms are also introduced in order to take the advantages of both $L^1$ and $L^2$ regularization.
For instance, \cite{taoetal2016} adopts the $l^{2,p}$ matrix norm as the regularization.
\par
In the following part, we state and prove the convergence results under the mild conditions (A.1)-(A.6).

\subsection{Convergence in Expectation}
\begin{theorem}
Let $(\mathbf{x}_t)_{t\ge0}$ be computed by Algorithm 1.
Suppose that the conditions (A.1)-(A.5) hold and the stepsize satisfies (A.6).
Then
\begin{equation*}
\sum_{t=0}^\infty \gamma_t {\mathbb{E} \Vert \nabla f(\mathbf{x}_t) \Vert}^2 \le C_0 ,
\end{equation*}
where $C_0$ is a constant depending on $f, \mathbf{x}_0, \beta, s, (\gamma_t), C, \sigma$.
\end{theorem}
Instead of working over $(\mathbf{x}_t)_{t \ge 0}$ directly, which can be complicated and hinder the intuitions, we utilize immediate variables to simplify the presentation and facilitate the analysis.
We start by defining $(\mathbf{p}_t)_{t\ge0}$ and deriving recursive formulas.
\begin{lemma}
Define
\begin{equation*}
\textbf{p}_t \ \hat{=} \ \begin{cases}
0, & t = 0;\\
\dfrac{\beta}{1\!-\!\beta}(\mathbf{x}_t - \mathbf{x}_{t-1} + s \gamma_{t-1} \mathbf{g}_{t-1}), & t \ge 1,
\end{cases}
\end{equation*}
and
\begin{equation*}
\mathbf{v}_t \ \hat{=} \ \frac{1\!-\!\beta}{\beta} \mathbf{p}_t, \quad t \ge 0 .
\end{equation*}
Let $\mathbf{z}_t = \mathbf{x}_t + \mathbf{p}_t$.
Then for any $t \ge 0$, we have
\begin{align*}
\mathbf{z}_{t+1} &= \mathbf{z}_t - \frac{\gamma_t}{1\!-\!\beta} \mathbf{g}_t,\\
\mathbf{v}_{t+1} &= \beta \mathbf{v}_t + ((1\!-\!\beta)s-1) \gamma_t \mathbf{g}_t.
\end{align*}
\end{lemma}
\begin{remark}
The proof of the recursion is straightforward and can be found in \cite{yanetal2018}, although only constant stepsizes are processed there.
So we omit this proof.
A similar recursion is given by \cite{ghadimietal2015} with $s = 0$ or $1$.
Note that the lemma does not hold when $\beta$ is replaced by $(\beta_t)$.
In the case of SHB with a variable momentum factor $(\beta_t)$, a recursion can be found in \cite{chenetal2019}.
\end{remark}
\begin{remark}
Out of aesthetics and succinctness, we use the SUM framework established by Yan et al. in \cite{yanetal2018} to obtain the convergence of SGD, SNAG, and SHB, but the unification is not essential to our demonstration.
Actually, we can discuss these stochastic gradient algorithms one by one along the same proof line.
\end{remark}
Applying Newton-Leibniz formula, we have the following estimate of a $C^{1+\alpha}$ function.
\begin{lemma}
Let $f \in C^{1+\alpha}(\mathbb{R}^d)$ with $A > 0$ being the H{\"o}lder index of the gradient $\nabla f$.
Then for any $\mathbf{x}, \mathbf{z} \in \mathbb{R}^d$,
\begin{equation*}
f(\mathbf{x}+\mathbf{z}) - f(\mathbf{x}) \le \langle\mathbf{z},\nabla f(\mathbf{x})\rangle + A \frac{\;{\Vert\mathbf{z}\Vert}^{1+\alpha}}{1\!+\!\alpha} .
\end{equation*}
\end{lemma}
Having the above preliminaries, we prove the convergence of the gradient in $L^2 (\Omega)$ space.
Applying common techniques to the $\alpha$-H{\"o}lder case, this argumentation bears a resemblance to most discussions of stochastic gradient methods (\cite{ghadimilan2013}, \cite{ghadimilan2016}, \cite{yanetal2018}, \cite{cutkoskyorabona2019}, etc.), where the rate of convergence of the form $\min\mathbb{E}{\Vert\nabla f(\mathbf{x}_t)\Vert}^2$ or $\mathbb{E}{\Vert\nabla f(\mathbf{x}_{R_t})\Vert}^2$ is obtained by an approximation of $\mathbb{E}[f(\mathbf{x}_{t+1})-f(\mathbf{x}_0)]$.
\par\vspace{2pt}
\begin{proof}[Proof of Theorem 1]
Define $\mathbf{\delta}_t \ \hat{=} \ \mathbf{g}_t - \nabla f(\mathbf{x}_t),\ t \ge 0$.
By Jessen's inequality, for any $a, b \ge 0$, $(a + b)^{1+\alpha} \le 2^\alpha (a^{1+\alpha} + b^{1+\alpha}) \le 2 (a^{1+\alpha} + b^{1+\alpha})$.
According to Lemma 1, noting the recursion of $(\mathbf{z}_t)$, we have, for any $t \in \mathbb{N}$,
\begin{align}
&f(\mathbf{z}_{t+1}) - f(\mathbf{z}_t) \nonumber\\
\le &- \frac{\gamma_t}{1\!-\!\beta} \langle\nabla f(\mathbf{z}_t), \mathbf{g}_t\rangle + \frac{A \gamma_t^{1+\alpha}}{(1\!+\!\alpha)(1\!-\!\beta)^{1+\alpha}} {\Vert\mathbf{g}_t\Vert}^{1+\alpha} \nonumber\\
\le &- \frac{\gamma_t}{1\!-\!\beta} \langle\nabla f(\mathbf{z}_t), \mathbf{g}_t\rangle + \frac{2A\gamma_t^{1+\alpha}}{(1\!+\!\alpha)(1\!-\!\beta)^{1+\alpha}} ({\Vert\nabla f(\mathbf{x}_t)\Vert}^{1+\alpha} \!+\! {\Vert\mathbf{\delta}_t\Vert}^{1+\alpha}).
\end{align}
Taking conditional expectations with respect to $\mathcal{F}_t$ on both sides, under assumptions (A.3)-(A.5), we obtain
\begin{align}
&\mathbb{E}[f(\mathbf{z}_{t+1}) - f(\mathbf{z}_t) |\mathcal{F}_t] \nonumber\\
\le &- \frac{\gamma_t}{1\!-\!\beta} \langle\nabla f(\mathbf{z}_t), \nabla f(\mathbf{x}_t)\rangle + \frac{2A\gamma_t^{1+\alpha}}{(1\!+\!\alpha)(1\!-\!\beta)^{1+\alpha}} ({\Vert\nabla f(\mathbf{x}_t)\Vert}^{1+\alpha} \!+\! \mathbb{E}[{\Vert\mathbf{\delta}_t\Vert}^{1+\alpha}|\mathcal{F}_t]) \nonumber\\
\le &-\frac{\gamma_t}{1\!-\!\beta} \langle\nabla f(\mathbf{z}_t), \nabla f(\mathbf{x}_t)\rangle + \frac{2A\gamma_t^{1+\alpha}}{(1\!+\!\alpha)(1\!-\!\beta)^{1+\alpha}} ((C^{\frac{1+\alpha}{2}}\!+\!1){\Vert\nabla f(\mathbf{x}_t)\Vert}^{1+\alpha} \!+\! \sigma^{1+\alpha}) \nonumber\\
= &- \frac{\gamma_t}{1\!-\!\beta} \langle\nabla f(\mathbf{z}_t) \!-\! \nabla f(\mathbf{x}_t), \nabla f(\mathbf{x}_t)\rangle - \frac{\gamma_t}{1\!-\!\beta} {\Vert\nabla f(\mathbf{x}_t)\Vert}^2 \nonumber\\[-2pt]
&+ \frac{2A\gamma_t^{1+\alpha}}{(1\!+\!\alpha)(1\!-\!\beta)^{1+\alpha}} ((C^{\frac{1+\alpha}{2}}\!+\!1){\Vert\nabla f(\mathbf{x}_t)\Vert}^{1+\alpha} \!+\! \sigma^{1+\alpha}) \nonumber\\
\le &\frac{A\gamma_t^{1+\alpha}}{2(1\!-\!\beta)^{1+\alpha}} {\Vert\nabla f(\mathbf{x}_t)\Vert}^2 + \frac{\gamma_t^{1-\alpha}}{2A(1\!-\!\beta)^{1-\alpha}} {\Vert\nabla f(\mathbf{z}_t) \!-\! \nabla f(\mathbf{x}_t)\Vert}^2 - \frac{\gamma_t}{1\!-\!\beta} {\Vert\nabla f(\mathbf{x}_t)\Vert}^2 \nonumber\\[-2pt]
&+ \frac{2A\gamma_t^{1+\alpha}}{(1\!+\!\alpha)(1\!-\!\beta)^{1+\alpha}} ((C^{\frac{1+\alpha}{2}}\!+\!1){\Vert\nabla f(\mathbf{x}_t)\Vert}^{1+\alpha} \!+\! \sigma^{1+\alpha}),
\end{align}
where the second inequality follows from Jessen's inequality for conditional expectations and the last inequality follows from Cauchy–Schwarz inequality.
\par
Define $\Gamma_t \ \hat{=} \ \sum_{i=0}^{t}\beta^i = \frac{1\!-\!\beta^{t+1}}{1\!-\!\beta}$.
Observe that
\begin{align}
&{\Vert\nabla f(\mathbf{z}_t) \!-\! \nabla f(\mathbf{x}_t)\Vert}^2 \nonumber\\
\le &A^2{\Vert\mathbf{p}_t\Vert}^{2\alpha} \nonumber\\
= &\frac{A^2\beta^{2\alpha}}{(1\!-\!\beta)^{2\alpha}} {\Vert\mathbf{v}_t\Vert}^{2\alpha} \nonumber\\
= &\frac{A^2\beta^{2\alpha}{|(1\!-\!\beta)s-1|}^{2\alpha}}{(1\!-\!\beta)^{2\alpha}} {\left\| \sum_{i=0}^{t-1} \beta^i \, \gamma_{t-1-i} \, \mathbf{g}_{t-1-i} \right\|}^{2\alpha} \nonumber\\
\le &\frac{A^2\beta^{2\alpha}{|(1\!-\!\beta)s-1|}^{2\alpha}}{(1\!-\!\beta)^{2\alpha}} \Gamma_{t-1}^\alpha \left\{ \sum_{i=0}^{t-1} \beta^i \, \gamma_{t-1-i}^2 {\Vert\mathbf{g}_{t-1-i}\Vert}^2 \right\}^\alpha \nonumber\\
\le &\frac{A^2\beta^{2\alpha}{|(1\!-\!\beta)s-1|}^{2\alpha}}{(1\!-\!\beta)^{2\alpha}} \Gamma_{t-1}^\alpha \sum_{i=0}^{t-1} \beta^{i\alpha} \, \gamma_{t-1-i}^{2\alpha} \, {\Vert\mathbf{g}_{t-1-i}\Vert}^{2\alpha} \nonumber\\
\le &\frac{A^2\beta^{2\alpha}{|(1\!-\!\beta)s-1|}^{2\alpha}}{(1\!-\!\beta)^{3\alpha}} \sum_{i=0}^{t-1} \beta^{i\alpha} \, \gamma_{t-1-i}^{2\alpha} \, {\Vert\mathbf{g}_{t-1-i}\Vert}^{2\alpha},
\end{align}
where the second equation follows from the recursion of $(\mathbf{v}_t)$ and the second inequality follows from H{\"o}lder's inequality for probability measures.
\par
Now we substitute the second term of RHS of (9) with (10) and obtain
\begin{align}
&\mathbb{E}[f(\mathbf{z}_{t+1}) - f(\mathbf{z}_t)] \nonumber\\
\le &\frac{A\beta^{2\alpha}{|(1\!-\!\beta)s-1|}^{2\alpha}}{(1\!-\!\beta)^{1+2\alpha}} \Bigg\{ \sum_{i=0}^{t-1} {\beta}^{i\alpha} \,\gamma_t^{1-\alpha} \, \gamma_{t-1-i}^{2\alpha} [\sigma^{2\alpha} + (C^\alpha\!+\!1)\mathbb{E}{\Vert\nabla f(\mathbf{x}_{t-1-i})\Vert}^{2\alpha}] \Bigg\} \nonumber\\
&+ \left( \frac{A\gamma_t^{1+\alpha}}{2(1\!-\!\beta)^{1+\alpha}} - \frac{\gamma_t}{1\!-\!\beta} \right) \mathbb{E}{\Vert\nabla f(\mathbf{x}_t)\Vert}^2 + \frac{2A\gamma_t^{1+\alpha}}{(1\!+\!\alpha)(1\!-\!\beta)^{1+\alpha}} (C^{\frac{1+\alpha}{2}}\!+\!1) \mathbb{E}{\Vert\nabla f(\mathbf{x}_t)\Vert}^{1+\alpha} \nonumber\\
&+ \frac{2A\sigma^{1+\alpha}\gamma_t^{1+\alpha}}{(1\!+\!\alpha)(1\!-\!\beta)^{1+\alpha}}.
\end{align}
\par
Our next step is to bound the summation of the first terms of RHS of (11) as below.
Let $S$ denote $1+ \sum_{t=0}^\infty \gamma_t^{1+\alpha}$ and we have
\begin{align}
&\sum_{t=0}^\infty \sum_{i=0}^{t-1} \beta^{i\alpha} \gamma_t^{1-\alpha} \gamma_{t-1-i}^{2\alpha} \mathbb{E}{\Vert\nabla f(\mathbf{x}_{t-1-i})\Vert}^{2\alpha} \nonumber\\
= &\sum_{i=0}^\infty \beta^{i\alpha} \sum_{t=i+1}^\infty \gamma_t^{1-\alpha} \Big( \gamma_{t-1-i}^{2\alpha} \mathbb{E}{\Vert\nabla f(\mathbf{x}_{t-1-i})\Vert}^{2\alpha} \Big) \nonumber\\
\le &\sum_{i=0}^\infty \beta^{i\alpha} \Bigg\{ \Big( \sum_{t=i+1}^\infty \gamma_t^{1+\alpha} \Big)^\frac{1-\alpha}{1+\alpha} \Big( \sum_{t=i+1}^\infty \gamma_{t-1-i}^{1+\alpha} \mathbb{E}{\Vert\nabla f(\mathbf{x}_{t-1-i})\Vert}^{1+\alpha} \Big)^\frac{2\alpha}{1+\alpha} \Bigg\} \nonumber\\
\le &\sum_{i=0}^\infty \beta^{i\alpha} \Big( 1 \!+\! \sum_{t=0}^\infty \gamma_t^{1+\alpha} \Big) \Big( 1 \!+\! \sum_{t=0}^\infty \gamma_t^{1+\alpha}\mathbb{E}{\Vert\nabla f(\mathbf{x}_{t})\Vert}^{1+\alpha} \Big) \nonumber\\
\le &\frac{S}{1\!-\!\beta^\alpha} + \frac{S}{1\!-\!\beta^\alpha} \sum_{t=0}^\infty \gamma_t^{1+\alpha} \mathbb{E}{\Vert\nabla f(\mathbf{x}_{t})\Vert}^{1+\alpha},
\end{align}
where in the first equality we have used H{\"o}lder's inequality and Jessen's inequality, and in the second one ${|a|}^c \le 1 + |a|$ for any $c \in [0,1]$.
Similarly,
\begin{equation}
\sum_{t=0}^\infty \sum_{i=0}^{t-1} \beta^{i\alpha} \gamma_t^{1-\alpha} \gamma_{t-1-i}^{2\alpha} \le \frac{1}{1\!-\!\beta^\alpha} \sum_{t=0}^\infty \gamma_t^{1+\alpha} \le \frac{S}{1\!-\!\beta^\alpha} .
\end{equation}
\par
For simplicity, we set
\begin{align*}
&C_1(S) \ \hat{=} \ \frac{A}{2(1\!-\!\beta)^{1+\alpha}} + \frac{2A(C^{\frac{1+\alpha}{2}}\!+\!1)}{(1\!+\!\alpha)(1\!-\!\beta)^{1+\alpha}} + \frac{A\beta^{2\alpha}{|(1\!-\!\beta)s-1|}^{2\alpha}S(C^\alpha\!+\!1)}{(1\!-\!\beta)^{1+2\alpha}(1\!-\!\beta^\alpha)} ,\\
&C_2 \ \hat{=} \ \frac{A\beta^{2\alpha}{|(1\!-\!\beta)s-1|}^{2\alpha}(\sigma^{2\alpha}\!+\!2C^\alpha\!+\!2)}{(1\!-\!\beta)^{1+2\alpha}(1\!-\!\beta^\alpha)} + \frac{2A(C^{\frac{1+\alpha}{2}}\!+\!1)}{(1\!+\!\alpha)(1\!-\!\beta)^{1+\alpha}} + \frac{2A\sigma^{1+\alpha}}{(1\!+\!\alpha)(1\!-\!\beta)^{1+\alpha}} .
\end{align*}
\par
Then, we can sum up the inequalities involving $\mathbb{E} [f(\mathbf{z}_{t+1})-f(\mathbf{z}_t)]$ and rearrange the terms.
Noting that ${|a|}^c \le 1 + {|a|}^2$ for any $c \in [0,2]$, we immediately obtain the following succinct formula.
\begin{equation}
\sum_{t=0}^\infty \Big( \frac{\gamma_t}{1\!-\!\beta} - C_1(S)\gamma_t^{1+\alpha} \Big) \mathbb{E}\Vert\nabla f(\mathbf{x}_t)\Vert^2 \le (f(\mathbf{x}_0) -f_*) +SC_2 .
\end{equation}
After some finite number of terms, the $t$-th term of LHS is bigger than $\frac{\gamma_{t-1}}{2(1\!-\!\beta)}\mathbb{E}{\Vert\nabla f(\mathbf{x}_{t-1})\Vert}^2$.
Moreover, note $\mathbb{E}{\Vert\nabla f(\mathbf{x}_t)\Vert}^2 < \infty$ for any $t$, which can be shown by induction.
We conclude
\begin{equation}
\sum_{t=0}^\infty \frac{\gamma_t}{2(1\!-\!\beta)} \mathbb{E}{\Vert\nabla f(\mathbf{x}_t)\Vert}^2 \le C' .
\end{equation}
\end{proof}
\begin{corollary}
Let $R_t$ be a random variable taking on a value in $\{0, 1, \dots, t\}$ with the probability measure
\begin{equation*}
\mathbb{P}(R_t=i) = \frac{\gamma_i}{\sum_{i=0}^t \gamma_i}, \quad 0 \le i \le t .
\end{equation*}
And $(R_t)$ is independent of $(\xi_t)$.
Then, under assumptions (A.1)-(A.6),
\begin{equation*}
\mathbb{E}{\Vert\nabla f(\mathbf{x}_{R_t})\Vert}^2 \le \frac{C_0}{\sum_{i=0}^t \gamma_i} \to 0 \quad \textnormal{as} \ t \to 0 .
\end{equation*}
Furthermore, if more information on the stepsize is given, we can obtain an explicit upper bound on the rate of convergence.
\end{corollary}
\begin{corollary}
The estimates still hold under assumptions (A.1)-(A.5) for a finite number of iterations.
With a given positive integer $t$, consider $(\mathbf{x}_i)_{i=0,1,\dots,t}$ which is computed by Algorithm 1 with $\gamma_i = 1/t^{\frac{1}{1+\alpha}}$ for $i = 0,\dots,t\!-\!1$.
Noting $\sum_{i=0}^{t-1} \gamma_i = t^\frac{\alpha}{1+\alpha}$ and $\sum_{i=0}^{t-1} (\gamma_i)^{1+\alpha} = 1$, we immediately have $S=2$ and
\begin{equation*}
\sum_{i=0}^{t-1} \Big( \frac{\gamma_i}{1\!-\!\beta} - C_1(2)\gamma_i^{1+\alpha} \Big) \mathbb{E}\Vert\nabla f(\mathbf{x}_i)\Vert^2 \le (f(\mathbf{x}_0) -f_*) +2C_2 .
\end{equation*}
Following the notation $R_{t-1}$ of Corollary 3.1, when $t \ge (2(1\!-\!\beta)C_1(2))^\frac{1+\alpha}{\alpha}$, we have that
\begin{equation*}
\mathbb{E}\Vert\nabla f(\mathbf{x}_{R_{t-1}})\Vert^2 \le \frac{2(1\!-\!\beta)((f(\mathbf{x}_0) \!-\!f_*) \!+\!2C_2)}{\sum_{i=0}^{t-1} \gamma_i} = C_0 t^{-\frac{\alpha}{1+\alpha}} .
\end{equation*}
The above inequality shows that momentum methods ensure $\mathbb{E}{\Vert\nabla f\Vert}^2 \le \epsilon$ in $\mathcal{O}(1/\epsilon^{\frac{1+\alpha}{\alpha}})$ iterations with the constant stepsize being proportional to $\epsilon^\frac{1}{\alpha}$, which is aligned with our intuition.
\end{corollary}

\subsection{Almost Sure Convergence}
Based on the result in $L^2$ sense, almost sure convergence can be established.
\begin{theorem}
Let $(\mathbf{x}_t)_{t\ge0}$ be computed by Algorithm 1.
Also suppose that conditions (A.1)-(A.5) are fulfilled and stepsizes are chosen such that (A.6) holds.
Then, we have the following three assertions:
\vspace{-4pt}
\begin{enumerate}[(i)]
\item the sequence $f(\mathbf{x}_t)$ converges almost surely;
\vspace{-2pt}
\item the sequence $\nabla f(\mathbf{x}_t)$ converges to zero almost surely;
\vspace{-2pt}
\item if, in addition, $f$ has finite critical points in $\{ \mathbf{x} \in \mathbb{R}^d | a \le f(\mathbf{x}) \le b \}$ for any $a < b$ and $\lim_{\Vert\mathbf{x}\Vert \to \infty} f(\mathbf{x}) = +\infty$, then, $\mathbf{x_t}$ converges almost surely.
\end{enumerate}
\end{theorem}
Theorem 2 is not a direct inference that follows from Theorem 1.
We actually need to carefully restrict the dynamics introduced by noise and momentum.
Then we can illustrate that $f(\mathbf{x}_t)$ is approximately a supermartingale and construct subsequent arguments.
\par
The following lemma is of great importance in our argumentation, since it bridges the gap between $\mathbf{z}_t$ and $\mathbf{x}_t$ by demonstrating that $\mathbf{p}_t$ converges to $0$.
To some extent, this original observation enables us to treat momentum algorithms as vanilla ones without momentum.
\begin{lemma}
With the notation before, under assumptions (A.1)-(A.6), we have
\vspace{-4pt}
\begin{enumerate}[(i)]
\item $\sum_{t=0}^\infty \gamma_t {\Vert\mathbf{p}_t\Vert}^{2\alpha} < \infty \,$ a.s.;
\vspace{-2pt}
\item the sequence $\mathbf{p}_t$ converges to zero almost surely.
\end{enumerate}
\end{lemma}
\begin{proof}
As in (10), we	use the same approximation and obtain
\begin{align}
&\mathbb{E}{\Vert\mathbf{p}_t\Vert}^{2\alpha} \nonumber\\
\le &\frac{\beta^{2\alpha}{|(1\!-\!\beta)s-1|}^{2\alpha}}{(1\!-\!\beta)^{3\alpha}} \sum_{i=0}^{t-1} \beta^{i\alpha} \, \gamma_{t-1-i}^{2\alpha} \, \mathbb{E}{\Vert\mathbf{g}_{t-1-i}\Vert}^{2\alpha} \nonumber\\
\le &\frac{\beta^{2\alpha}{|(1\!-\!\beta)s-1|}^{2\alpha}(C^\alpha\!+\!1)}{(1\!-\!\beta)^{3\alpha}} \Bigg\{ \sum_{i=0}^{t-1} \beta^{i\alpha} \, \gamma_{t-1-i}^{2\alpha} \mathbb{E}{\Vert\nabla f(\mathbf{x}_{t-1-i})\Vert}^{2\alpha} \Bigg\} \nonumber\\
&+ \frac{\beta^{2\alpha}{|(1\!-\!\beta)s-1|}^{2\alpha}\sigma^{2\alpha}}{(1\!-\!\beta)^{3\alpha}} \sum_{i=0}^{t-1} \beta^{i\alpha} \, \gamma_{t-1-i}^{2\alpha}.
\end{align}
\par
We use $C'$ to denote a constant.
Multiplying both sides of (16) by $\gamma_t$ and summarizing over t, we have
\begin{align}
&\mathbb{E} \left[ \sum_{t=0}^\infty \gamma_t {\Vert\mathbf{p}_t\Vert}^{2\alpha} \right] \nonumber\\
\le &C'\sum_{t=0}^\infty \gamma_t \sum_{i=0}^{t-1} \beta^{i\alpha} \, \gamma_{t-1-i}^{2\alpha} \, \mathbb{E}{\Vert\nabla f(\mathbf{x}_{t-1-i})\Vert}^{2\alpha} + C'\sum_{t=0}^\infty \gamma_t \sum_{i=0}^{t-1} \beta^{i\alpha} \, \gamma_{t-1-i}^{2\alpha} \nonumber\\
\le &C'\sum_{i=0}^\infty \beta^{i\alpha} \Big( \sum_{t=i+1}^\infty \gamma_t^\frac{1+\alpha}{1-\alpha} \Big)^\frac{1-\alpha}{1+\alpha} \Big( \sum_{t=i+1}^\infty \gamma_{t-1-i}^{1+\alpha} \mathbb{E}{\Vert\nabla f(\mathbf{x}_{t-1-i})\Vert}^{1+\alpha} \Big)^\frac{2\alpha}{1+\alpha} \nonumber\\
&+ C'\sum_{i=0}^\infty \beta^{i\alpha} \Big( \sum_{t=i+1}^\infty \gamma_t^\frac{1+\alpha}{1-\alpha} \Big)^\frac{1-\alpha}{1+\alpha} \Big( \sum_{t=i+1}^\infty \gamma_{t-1-i}^{1+\alpha}\Big)^\frac{2\alpha}{1+\alpha} \nonumber\\
< &\infty,
\end{align}
which is a direct application of Theorem 1 and (A.6).
Therefore, $\sum_{t=0}^\infty \gamma_t {\Vert\mathbf{p}_t\Vert}^{2\alpha} < \infty$ a.s.
\par
Similar to the derivation of (i), we have
\begin{align}
&\mathbb{E} \left[ \sum_{t=0}^\infty {\Vert\mathbf{p}_t\Vert}^2 \right] \nonumber\\
\le &C' \sum_{i=0}^\infty \beta^i \sum_{t=i+1}^\infty \gamma_{t-1-i}^2 \mathbb{E}{\Vert\nabla f(\mathbf{x}_{t-1-i})\Vert}^2 + C' \sum_{i=0}^\infty \beta^i \sum_{t=i+1}^\infty \gamma_{t-1-i}^2 \nonumber\\
< &\infty.
\end{align}
We conclude that $\sum_{t=0}^\infty {\Vert\mathbf{p}_t\Vert}^2 < \infty$ and $\Vert\mathbf{p}_t\Vert \to 0$ almost surely.
\end{proof}
\par
Now, we are ready to prove the almost sure convergence of the SUM method.
\begin{proof}[Proof of Theorem 2]
Noting $\Vert\mathbf{x}_{t+1} \!-\! \mathbf{x}_t\Vert \to 0$, we only need to prove assertion (i) and assertion (ii), because assertion (iii) is a direct consequence.
The proof consists of three steps whose main body is along the lines of the proof of Proposition 4.1 in \cite{bertsekastsitsiklis1996}.
First we use Doob's supermartingale convergence theorem to bound the effects of the noise and obtain the result that $f(\mathbf{z}_t)$ converges and $\liminf\limits \, \Vert\nabla f(\mathbf{x}_t)\Vert = 0$ a.s.
Assuming $\limsup\limits \, \Vert\nabla f(\mathbf{x}_t)\Vert > 0$, we then proceed to a detailed discussion of upcrossing intervals and reduce this assumption to absurdity.
The remaining part, i.e., $f(\mathbf{x}_t)$ converges almost surely, is completed in the last step.
\par
In truth, we can take the first step swiftly by means of the Robbins-Siegmund lemma, but we provide a self-contained derivation here.
\par
As is similar to the approximation of $\mathbb{E}[f(\mathbf{z}_{t+1}) \!-\!  f(\mathbf{z}_t) |\mathcal{F}_t]$ in (9) and (10), we have
\begin{equation}
\mathbb{E}[f(\mathbf{z}_{t+1})|\mathcal{F}_t] \le f(\mathbf{z}_t) -X_t +Y_t +Z_t ,
\end{equation}
with
\begin{align*}
X_t \ &\hat{=} \ \left( \frac{\gamma_t}{2(1\!-\!\beta)} -\frac{2A(1\!+\!C)\gamma_t^{1+\alpha}}{(1\!+\!\alpha)(1\!-\!\beta)^{1+\alpha}} \right) {\Vert\nabla f(\mathbf{x}_t)\Vert}^2 , \\
Y_t \ &\hat{=} \ \frac{A^2\gamma_t}{2(1\!-\!\beta)} {\Vert\mathbf{p}_t\Vert}^{2\alpha} , \\
Z_t \ &\hat{=} \ \frac{A(2\sigma^2\!+\!1)}{(1\!+\!\alpha)(1\!-\!\beta)^{1+\alpha}} \gamma_t^{1+\alpha} .
\end{align*}
Note that $(Y_t)$ and $(Z_t)$ are nonnegative $(\mathcal{F}_t)$-adapted processes.
Also, $(X_t)$ is adapted and there exists a constant $T_X$ such that $X_t \ge 0$ for any $t > T_X$.
Consider the adapted process $(\widetilde{f}_t)$ defined by $\widetilde{f}_t \ \hat{=} \ f(\mathbf{z}_t) +\sum_{i=0}^{t-1}X_i -\sum_{i=0}^{t-1}Y_i -\sum_{i=0}^{t-1}Z_i $, for any $t \ge 0$.
The above inequality can be written as
\begin{equation}
\mathbb{E}[\widetilde{f}_{t+1}|\mathcal{F}_t] \le \widetilde{f}_t ,
\end{equation}
which means that $(\widetilde{f}_t)$ is a $(\mathcal{F}_t)$-supermartingale.
\par
In fact, the expectation of negative part of $(\widetilde{f}_t)$ is bounded.
We make use of (A.6) and (17) and obtain
\begin{align}
&\mathbb{E}[(\widetilde{f}_t)^-] \nonumber\\
\le &(f(\mathbf{z}_t))^- +\sum_{i=0}^{T_X}\mathbb{E}[(X_i)^-] +\sum_{i=0}^{t-1}\mathbb{E}[Y_i] +\sum_{i=0}^{t-1}Z_i \nonumber\\
= &(f(\mathbf{z}_t))^- +\sum_{i=0}^{T_X}\mathbb{E}[(X_i)^-] +\frac{A^2}{2(1\!-\!\beta)}\sum_{i=0}^{t-1}\gamma_i\mathbb{E}{\Vert\mathbf{p}_i\Vert}^{2\alpha} +\frac{A(2\sigma^2\!+\!1)}{(1\!+\!\alpha)(1\!-\!\beta)^{1+\alpha}}\sum_{i=0}^{t-1}\gamma_i^{1+\alpha} \nonumber\\
\le &C',
\end{align}
where the constant $C'$ does not depend on $t$.
\par
Subsequently, using Doob's supermartingale convergence theorem, we deduce that there exists a random variable $\eta \in L^1(\Omega, \mathcal{F}, \mathbb{P})$ such that $\widetilde{f}_t \to \eta$ a.s.
Thanks to the convergence of the sum of $Y_t$ and $Z_t$, $(f(\mathbf{z}_t) +\sum_{i=0}^{t-1}X_i)$ converges almost surely.
Note that $\inf_{\mathbf{x}}f(\mathbf{x}) > -\infty$ and $(\sum_{i=0}^{t-1}X_i)$ is entirely non-decreasing after $T_X$ terms.
As a result,
\begin{equation*}
\lim_{t\to\infty} f(\mathbf{z}_t) \ \textnormal{exists and} \ \sum_{t=0}^\infty X_t < \infty \ \textnormal{a.s.}
\end{equation*}
Since $\sum_{t=0}^\infty \gamma_t = \infty$, we immediately obtain
\begin{equation*}
\liminf\limits_{t\to\infty} \, \Vert\nabla f(\mathbf{x}_t)\Vert = 0 \ \textnormal{a.s.}
\end{equation*}
The first step has been accomplished.
\par
\vspace{2pt}
We say that the time interval $\{t, t+1, \dots, \overline{t}\}$ is an \emph{upcrossing interval} of $(\Vert\nabla f(\mathbf{x}_s)\Vert)$ from $a$ to $b$, if $\Vert\nabla f(\mathbf{x}_t)\Vert < a$, $\Vert\nabla f(\mathbf{x}_{\overline{t}})\Vert > b$, and $a \le \Vert\nabla f(\mathbf{x}_s)\Vert \le b$ for $t<s<\overline{t}$.
Now we assume the opposite of the proposition that $\nabla f(\mathbf{x}_t)$ converges to zero a.s. and show its nonsense.
\par
Assume there exist a constant $\epsilon > 0$ and a measurable set $\Omega_1$ with $\mathbb{P}(\Omega_1) > 0$ such that, for any $\omega \in \Omega_1$, $\limsup\limits \, \Vert\nabla f(\mathbf{x}_t)\Vert (\omega) \ge 2\epsilon$.
The fact that the limit inferior of $\nabla f(\mathbf{x}_t)$ is zero and the above assumption lead to an infinite number of upcrossings from $\frac{\epsilon}{2}$ to $\epsilon$ for $\omega \in \Omega_1$.
Denote the $k$-th upcrossing interval by $\{t_k, t_{k+1}, \dots, \overline{t}_k\}$.
\par
Define
\begin{equation*}
\chi_t \ \hat{=} \ \mathds{1}_{\{ \Vert\nabla f(\mathbf{x}_t)\Vert \le \epsilon \}} \in \mathcal{F}_t ,
\end{equation*}
and
\begin{equation*}
\mathbf{u}_t \ \hat{=} \ \sum_{i=0}^{t-1} \chi_i \gamma_i (\mathbf{g}_i \!-\! \mathbb{E}[\mathbf{g}_i|\mathcal{F}_i]) = \sum_{i=0}^{t-1} \chi_i\gamma_i\mathbf{\delta}_i \in \mathcal{F}_t .
\end{equation*}
Then we proceed analogously as in the analysis of $(\widetilde{f}_t)$ and obtain the almost sure convergence of the $(\mathcal{F}_t)$-martingale $(\mathbf{u}_t)$.
Since $(\chi_t)$ is adapted and the stochastic oracle is unbiased, we have, for any $t \ge 0$,
\par
\begin{equation}
\mathbb{E}[\mathbf{u}_{t+1}|\mathcal{F}_t] = \mathbf{u}_t +\mathbb{E}[\chi_t\gamma_t\mathbf{\delta}_t|\mathcal{F}_t] = \mathbf{u}_t +\chi_t\gamma_t\mathbb{E}[\mathbf{\delta}_t|\mathcal{F}_t] =\mathbf{u}_t .
\end{equation}
It is also straightforward to verify the $L^2(\Omega)$ boundedness of $(\mathbf{u}_t)$.
\begin{align}
&\mathbb{E}{\Vert\mathbf{u}_{t+1}\Vert}^2 \nonumber\\
= &\mathbb{E}{\Vert\mathbf{u}_t\Vert}^2 + \mathbb{E}[ \mathbb{E}[{\Vert\chi_t\gamma_t\mathbf{\delta}_t\Vert}^2|\mathcal{F}_t] ] \nonumber\\
\le &\mathbb{E}{\Vert\mathbf{u}_t\Vert}^2 + \gamma_t^2\,\mathbb{E}[\chi_t\,(C{\Vert\nabla f(\mathbf{x}_t)\Vert}^2 \!+\! \sigma^2)] \nonumber\\
\le &\sum_{i=0}^t \gamma_i^2\,(C\epsilon^2 \!+\! \sigma^2).
\end{align}
Applying Doob's martingale convergence theorem to $(\mathbf{u}_t)$, we conclude that $\mathbf{u}_t$ converges a.s.
We immediately obtain the two limits:
\begin{align}
&\lim_{k\to\infty} \sum_{t=t_k}^{\overline{t}_k-1} \gamma_t\mathbf{\delta}_t = 0 \ \textnormal{for almost all}\ \omega \in \Omega_1 ;\\
&\lim_{k\to\infty} \gamma_{t_k}\mathbf{\delta}_{t_k} = 0 \ \textnormal{for almost all}\ \omega \in \Omega_1 .
\end{align}
\par
Using the recursion in Lemma 1, we have, for any $k > 0$,
\begin{align}
&\Vert\nabla f(\mathbf{x}_{t_k+1})\Vert - \Vert\nabla f(\mathbf{x}_{t_k})\Vert \nonumber\\
\le &\Vert\nabla f(\mathbf{z}_{t_k+1})\Vert - \Vert\nabla f(\mathbf{z}_{t_k})\Vert + A{\Vert\mathbf{p}_{t_k+1}\Vert}^\alpha + A{\Vert\mathbf{p}_{t_k}\Vert}^\alpha \nonumber\\
\le &A{\left\| \frac{\gamma_{t_k}}{1\!-\!\beta}\mathbf{g}_{t_k} \right\|}^\alpha + A{\Vert\mathbf{p}_{t_k+1}\Vert}^\alpha + A{\Vert\mathbf{p}_{t_k}\Vert}^\alpha \nonumber\\
\le &\frac{A\gamma_{t_k}^\alpha}{(1\!-\!\beta)^\alpha}{\Vert\nabla f(\mathbf{x}_{t_k})\Vert}^\alpha + \frac{A}{(1\!-\!\beta)^\alpha}{\Vert\gamma_{t_k}\mathbf{\delta}_{t_k}\Vert}^\alpha + A{\Vert\mathbf{p}_{t_k+1}\Vert}^\alpha + A{\Vert\mathbf{p}_{t_k}\Vert}^\alpha.
\end{align}
Note that the four terms on the right side tends to zero for almost all $\omega \in \Omega_1$, respectively.
Therefore, there exists a measurable set $\Omega_2$ such that $\Omega_2 \subset \Omega_1$, $\mathbb{P}(\Omega_2) = \mathbb{P}(\Omega_1) > 0$, and $\sum_{t=t_k}^{\overline{t}_k-1} \gamma_t\mathbf{\delta}_t (\omega) \to 0$, $\gamma_{t_k}\mathbf{\delta}_{t_k} (\omega) \to 0$, $\Vert\nabla f(\mathbf{x}_{t_k+1})\Vert(\omega) \!-\! \Vert\nabla f(\mathbf{x}_{t_k})\Vert(\omega) \to 0$, $\mathbf{p}_{t_k}(\omega) \to 0$, and $\mathbf{p}_{\overline{t}_k}(\omega) \to 0$ for any $\omega \in \Omega_2$.
\par
We arbitrarily fix $\omega \in \Omega_2$ and consider this sample path below.
Since above convergence properties, there exists a positive integer $K_1(\omega)$ such that, for any $k > K_1(\omega)$, $\Vert\nabla f(\mathbf{x}_{t_k})\Vert (\omega) \ge \frac{\epsilon}{4}$.
Additionally, by (A.2) and the fact both $\mathbf{p}_{t_k}(\omega)$ and $\mathbf{p}_{\overline{t}_k}(\omega)$ converge to zero, we can choose a large enough $K_2(\omega)$ such that, for $k > K_2(\omega)$,
\begin{align}
\frac{\epsilon}{4} &\le \Vert\nabla f(\mathbf{z}_{\overline{t}_k})\Vert(\omega) \!-\! \Vert\nabla f(\mathbf{z}_{t_k})\Vert(\omega) \nonumber\\
&\le A{\left\| \sum_{t=t_k}^{\overline{t}_k-1} \frac{\gamma_t}{1\!-\!\beta}\mathbf{g}_t \right\|}^\alpha \!(\omega) \nonumber\\
&\le \frac{A}{(1\!-\!\beta)^\alpha} {\left\|\sum_{t=t_k}^{\overline{t}_k-1}\gamma_t\mathbf{\delta}_t\right\|}^\alpha \!(\omega) \!+\! \frac{A}{(1\!-\!\beta)^\alpha} {\left\| \sum_{t=t_k}^{\overline{t}_k-1}\gamma_t\nabla f(\mathbf{x}_t) \right\|}^\alpha \!(\omega).
\end{align}
Since the first term on the right side tends to zero, we have
\begin{equation}
\liminf\limits_{k\to\infty} \sum_{t=t_k}^{\overline{t}_k-1} \gamma_t\Vert\nabla f(\mathbf{x}_t)\Vert(\omega) \ge \frac{(1\!-\!\beta)\epsilon^{1/\alpha}}{4^{1/\alpha}A^{1/\alpha}}.
\end{equation}
Noting that $\Vert\nabla f(\mathbf{x}_t)\Vert (\omega) \ge \frac{\epsilon}{4}$ for any $t \in [t_k(\omega),\, \overline{t}_k(\omega)-1]$ with $k > K_1(\omega)$, we obtain
\begin{equation}
\liminf\limits_{k\to\infty} \sum_{t=t_k}^{\overline{t}_k-1} \gamma_t{\Vert\nabla f(\mathbf{x}_t)\Vert}^2 (\omega) \ge \frac{(1\!-\!\beta)\epsilon^{1/\alpha}}{4^{1/\alpha}A^{1/\alpha}} \frac{\epsilon}{4} .
\end{equation}
\par
This immediately implies that $\sum_{t=0}^\infty \gamma_t {\Vert\nabla f(\mathbf{x}_t)\Vert}^2(\omega) = \infty$ for any $\omega \in \Omega_2$, which contradicts Theorem 1.
We conclude that $\nabla f(\mathbf{x}_t)$ converges to zero almost surely.
\par
\vspace{2pt}
Thus, it remains to show the convergence of $f(\mathbf{x}_t)$.
By Lemma 2, we have
\begin{equation}
|f(\mathbf{x}_t) - f(\mathbf{z}_t)| \le \Vert\nabla f(\mathbf{x}_t)\Vert \Vert\mathbf{p}_t\Vert + \frac{A}{1\!+\!\alpha}{\Vert\mathbf{p}_t\Vert}^{1+\alpha} \to 0 \ \textnormal{a.s.}
\end{equation}
Therefore, the convergence of $f(\mathbf{z}_t)$ implies that $f(\mathbf{x}_t)$ converges with probability 1.
\end{proof}
\begin{remark}
In fact, some conditions can be easily relaxed in our analysis.
Instead of strict restrictions on growth of $\nabla f$ and the unbiasedness of $\mathbf{g}$, we consider the following counterparts.
\par (A.2') $0 < \alpha \le 1$.
$\nabla f$ is $\alpha$-H{\"o}lder continuous and satisfies a linear growth condition, i.e., for some $A>0$,
\begin{equation*}
\Vert\nabla f(\mathbf{x})-\nabla f(\mathbf{y})\Vert \le A\,({\Vert\mathbf{x}-\mathbf{y}\Vert}^\alpha + {\Vert\mathbf{x}-\mathbf{y}\Vert}), \quad \forall \mathbf{x},\mathbf{y} \in \mathbb{R}^d.
\end{equation*}
\par (A.3') The conditional expectation of the stochastic oracle is not too short and makes an acute angle with the gradient of $f$.
More precisely, for some $K>0$,
\begin{equation*}
\langle\nabla f(\mathbf{x}_t),\,\mathbb{E}[\mathbf{g}_t|\mathcal{F}_t]\rangle \ge K{\Vert\nabla f(\mathbf{x}_t)\Vert}^2,\quad \forall t\ge0 .
\end{equation*}
\par (A.4') There exist positive constants $\sigma^2$ and $C$, such that
\begin{equation*}
\mathbb{E}[{\Vert\mathbf{g}_t\Vert}^2|\mathcal{F}_t] \le \sigma^2 + C {\Vert\nabla f(\mathbf{x}_t)\Vert}^2, \quad \forall t\ge0 .
\end{equation*}
\par
We can tune up our proof along the feasible route in this section.
For the sake of intuition, we just give the derivation under assumptions (A.1)-(A.6).
However, the whole key points of our proof do not change under assumptions (A.1), (A.2')-(A.4'), and (A.6).
Since these conditions still fit our analytical framework, we will get the same results.
\end{remark}

\section{Stochastic Stepsizes}
The convergence theorems in the previous section require deterministic stepsizes.
However, the performance fluctuates dramatically and heavily depends on the choice of stepsizes.
Towards the aim of obtaining easy-to-tune learning rates, adaptive algorithms such as AdaGrad, Adadelta, RMSprop, and Adam get employed and a considerable part of state-of-the-art results is achieved in deep learning articles.
See \cite{duchietal2011}, \cite{kingmaba2015}, \cite{reddietal2018}, \cite{chenetal2019}, and references therein.
So the determinacy of the stepsize appears to be restrictive and unnecessary from both practical and theoretical points of view.
\par
To fill the gap between deterministic and stochastic stepsizes, we will meticulously use an approach based upon a localization procedure.
The core of our discussion is to find a proper localizing sequence of stopping times that reduces $(\gamma_t)$ to a more regular one $(\gamma_t^{(N)})$ such that $\sum_{t=0}^\infty \gamma_t^{(N)} = \infty$ and $\sum_{t=0}^\infty (\gamma_t^{(N)})^{1+\alpha}$ $\le N$ a.s.
The remaining argumentation is produced alike.
Additionally, we obtain the almost sure convergence of some variants of adaptive algorithms with momentum as a by-product.
\par
Define $\mathcal{F}_t \ \hat{=} \ \sigma(\gamma_0, \xi_0, \gamma_1, \xi_1, \dots, \gamma_{t-1}, \xi_{t-1}, \gamma_t)$ for any $t\ge0$, i.e., $\mathcal{F}_t$ stands for the entire history of the algorithm up to and including the point at which $\gamma_t$ is selected, but before the update direction $\mathbf{g}(\mathbf{x}_t, \xi_t)$ is determined.
We have the following theorem.
\begin{theorem}
Let $(\mathbf{x}_t)_{t\ge0}$ be computed by Algorithm 1.
Suppose that the stepsize $(\gamma_t)$ is a nonnegative sequence satisfying $\sum_{t=0}^\infty \gamma_t =\infty$ and $\sum_{t=0}^\infty (\gamma_t)^{1+\alpha} <\infty$ a.s. and that $\xi_t$ is independent of $\mathcal{F}_t$ for any $t\ge0$.
Also, suppose the oracle $\mathbf{g}(\mathbf{x}_t, \xi_t)(\omega) = \mathbf{g}(\mathbf{x}_t(\omega), \xi_t(\omega))$, i.e., a sample of $\mathbf{g}$ is unaffected by the distribution of $\mathbf{x}$, and so is $\gamma$.
Then, under assumptions (A.1)-(A.5), the following hold:
\vspace{-4pt}
\begin{enumerate}[(i)]
\item the sequence $f(\mathbf{x}_t)$ converges almost surely;
\vspace{-2pt}
\item the sequence $\nabla f(\mathbf{x}_t)$ converges to zero almost surely;
\vspace{-2pt}
\item if, in addition, $f$ has finite critical points in $\{ \mathbf{x}\in\mathbb{R}^d | a \le f(\mathbf{x}) \le b\}$ for any $a < b$ and $\lim_{\Vert\mathbf{x}\Vert\to\infty} f(\mathbf{x}) = +\infty$, then, $\mathbf{x_t}$ converges almost surely.
\end{enumerate}
\end{theorem}
\begin{proof}
To establish the convergence we construct localizing stopping times to impose some additional assumptions about $(\gamma_t)$.
$C_1(\cdot)$ and $C_2$ are the same as those in the proof of Theorem 1.
Fix an integer $N \ge 0$.
Let $\gamma_*^{(N)} > 0$ denote the solution of the following equation:
\begin{equation}
\frac{\gamma}{2(1\!-\!\beta)} \ = \ C_1(N\!+\!1) \gamma^{1+\alpha} ,
\end{equation}
The left side of the above equation is greater than or equal to the right side for any $0\le\gamma\le\gamma_*^{(N)}$.
To begin with, we introduce the $(\mathcal{F}_t)$-stopping time:
\begin{equation*}
\tau_N \ \hat{=} \ \inf \Big\{ t\ge0: \sum_{i=0}^t \gamma_i^{1+\alpha} \ge N-\frac{\alpha\!+\!1}{\alpha} \Big\} .
\end{equation*}
Then, we adjust $(\gamma_t)$ to suit the needs of argumentation in the previous section.
Denote
\begin{equation*}
\gamma_t^{(N)} \ \hat{=} \ \gamma_t\mathds{1}_{\{t<\tau_N\}} +\frac{1}{t\!+\!1}\mathds{1}_{\{t\ge\tau_N\}} .
\end{equation*}
It is obvious that $(\gamma_t^{(N)})$ is adapted.
Now we fix an integer $M \ge 0$.
We subsequently define
\begin{equation*}
\gamma_t^{(N,M)} \ \hat{=} \ \begin{cases}
\gamma_t^{(N)}, & t < M;\\
\min \big(\gamma_t^{(N)},\,\gamma_*^{(N)}\big), & t \ge M.
\end{cases}
\end{equation*}
It is quite simple to verify that $(\gamma_t^{(N,M)})_{t\ge0}$ is adapted and satisfies that
\begin{equation}
\sum_{t=0}^\infty\gamma_t^{(N,M)}=\infty, \quad \sum_{t=0}^\infty(\gamma_t^{(N,M)})^{1+\alpha} \le N, \quad \textnormal{a.s.}
\end{equation}
\par
Therefore, by construction, the stepsize $(\gamma_t^{(N,M)})$ is endowed with the fine regularity.
Let $(\mathbf{x}_t^{(N,M)})$ be computed by Algorithm 1 with $(\gamma_t^{(N,M)})$ and $\mathbf{g}_t^{(N,M)}$ denote $\mathbf{g}(\mathbf{x}_t^{(N,M)}, \xi_t)$.
Define corresponding $\mathbf{p}_t^{(N,M)}$ in the same approach as in Lemma 1.
\par
Noting the independence of $\xi_t$ and $\mathcal{F}_t$, we argue as in the proof of Theorem 1 and obtain
\begin{align}
&\sum_{t=0}^\infty \mathbb{E} \Bigg[ \left( \frac{1}{1\!-\!\beta}\gamma_t^{(N,M)} - C_1(N\!+\!1) (\gamma_t^{(N,M)})^{1+\alpha} \right) {\Vert\nabla f(\mathbf{x}_t^{(N,M)})\Vert}^2 \Bigg] \nonumber\\
\le &(N\!+\!1)C_2 + (f(\mathbf{x}_0) - f_*).
\end{align}
Roughly speaking, $C_1(N\!+\!1)(\gamma_t^{(N.M)})^{1+\alpha}$ can be combined into $\frac{1}{1\!-\!\beta} \gamma_t^{(N.M)}$ because of the fact $0\le \gamma_t^{(N.M)} \le\gamma_*^{(N)}$ for any $t\ge M$.
More concretely, we have
\begin{equation}
\sum_{t=0}^\infty \mathbb{E}[\gamma_t^{(N,M)}{\Vert\nabla f(\mathbf{x}_t^{(N,M)})\Vert}^2] \le C(N,M).
\end{equation}
\par
Then we produce estimations and afterwards obtain the boundedness of $\mathbb{E} [\sum_{t=0}^\infty\gamma_t^{(N,M)}{\Vert\mathbf{p}_t^{(N,M)}\Vert}^{2\alpha}]$ and $\mathbb{E} [\sum_{t=0}^\infty{\Vert\mathbf{p}_t^{(N,M)}\Vert}^2]$, which echoes the derivation of (17) and (18).
As a result,
\begin{equation}
\sum_{t=0}^\infty \gamma_t^{(N,M)} {\Vert\mathbf{p}_t^{(N,M)}\Vert}^{2\alpha} < \infty \ \textnormal{a.s.}
\end{equation}
and
\begin{equation}
\mathbf{p}_t^{(N,M)} \to 0 \ \textnormal{a.s.}
\end{equation}
\par
Furthermore, we use the same approach as Theorem 2 based on the supermartingale convergence theorem and a sophisticated discussion of upcrossing intervals.
This part of the proof is unchanged.
Finally, we obtain the two limits:
\begin{align}
&\lim_{t\to\infty} \nabla f(\mathbf{x}_t^{(N,M)}) = 0 \ \textnormal{a.s.,} \\
&\lim_{t\to\infty} f(\mathbf{x}_t^{(N,M)}) \ \textnormal{exists} \ \textnormal{a.s.}
\end{align}
\par
Since $N$ and $M$ here are arbitrary, we are going to pass to the limit.
For an arbitrary $\omega$, except for some sample points forming a subset of a zero-probability event, we can choose a large enough integer $N(\omega)$ and then a proper $M(\omega)$ such that $(\gamma_t^{(N,M)})_{t\ge0}$ is identical to $(\gamma_t)_{t\ge0}$.
This immediately yields $(\mathbf{x}_t^{(N(\omega),M(\omega))}(\omega)) \equiv (\mathbf{x}_t(\omega))$ because a fixed sample of $\mathbf{g}$ and $\gamma$ is unaffected by other samples.
That is to say, the convergence of $f(\mathbf{x}_t)(\omega)$ and $\nabla f(\mathbf{x}_t)(\omega)$ holds.
As a matter of fact, we actually accomplish the proof.
\end{proof}
\par
\vspace{2pt}
Now we consider adaptive algorithms, i.e., $\gamma$ is a function of past stochastic gradients.
There is not enough research on this area, especially in the nonconvex setting.
What is worse, some algorithms have convergence issues.
For example, exponential moving average methods like Adam have flaws in an online convex setup, according to \cite{reddietal2018}.
\par
Thanks to Theorem 4.1, we can easily get mild sufficient conditions guaranteeing the convergence of some adaptive gradient algorithms.
We expect this perspective will shed a little light on adaptive stepsizes.
We consider the following generalized Adam algorithm as an example:
\begin{equation}
\left\{
\begin{aligned}
v_t &= \beta' v_{t-1} + (1\!-\!\beta'){\Vert\mathbf{g}_{t-1}\Vert}^2,\\
\mathbf{m}_{t+1} &= \beta\mathbf{m}_t - \frac{(1\!-\!\beta)}{(t\!+\!1)^{1\!/\!2+\epsilon}(\kappa\!+\!v_t^{1\!/\!2})}\mathbf{g}_t,\\
\mathbf{x}_{t+1} &= \mathbf{x}_t + \mathbf{m}_{t+1},
\end{aligned}
\right.
\end{equation}
where $v_0 = 0$, $\mathbf{m}_0 = 0$, $0\le\beta<1$, $0\le\beta'<1$, and $\kappa>0$ for ensuring the numerical stability of the stepsize.
Obviously, by setting $\gamma_t = (1\!-\!\beta)/((t\!+\!1)^{1\!/\!2+\epsilon}(\kappa\!+\!v_t^{1\!/\!2}))$ and $s=0$, the SUM method reduces to the above algorithm.
Note that the analysis here also applies to both a coordinate-wise stepsize and AdaFom (AdaGrad with First Order Momentum) firstly proposed by \cite{chenetal2019}.
A similar perturbed AdaGrad algorithm without momentum is analyzed in \cite{liorabona2018}.
\begin{corollary}
Consider the perturbed Adam (39).
Assume (A.1)-(A.3) and (A.5).
Suppose that $\frac{1}{1+\alpha}\!-\!\frac{1}{2} < \epsilon$ $\le \frac{1}{2}$ and there exists some constant $G>0$ such that $\Vert\nabla f(\mathbf{x}_t)\Vert \le G$ and $\Vert\mathbf{g}_t\Vert \le G$ for any $t$.
Then, the gradient converges to zero and the value of $f$ converges almost surely.
Moreover, $\liminf \, t^{1\!/\!2-\epsilon} {\Vert\nabla f(\mathbf{x}_t)\Vert}^2 = 0$ with probability 1.
\end{corollary}
\begin{proof}
It is clear to see that $\gamma_t \in \sigma(\xi_0, \dots, \xi_{t-1})$ for any $t$.
This yields $\xi_t$ is independent of $\mathcal{F}_t$.
It remains to show that effective stepsizes satisfy $\sum_{t=0}^\infty \gamma_t =\infty$ and $\sum_{t=0}^\infty (\gamma_t)^{1+\alpha} <\infty$ almost surely.
We only need to bound $\kappa\!+\!v_t^{1\!/\!2}$ as follows:
\begin{equation}
\kappa \le \kappa\!+\!v_t^{1\!/\!2} \le \kappa + \Big( \sum_{i=0}^{t-1}(\beta')^{t-1-i}(1\!-\!\beta')G^2 \Big)^{1\!/\!2} \le \kappa\!+\!G .
\end{equation}
\par
Thus, Theorem 3 can be utilized straightforwardly.
We know that $\sum_{t=0}^\infty\gamma_t{\Vert\nabla f(\mathbf{x}_t)\Vert}^2$ is finite with probability 1.
This means
\begin{equation}
\sum_{t=0}^\infty \frac{1}{(t\!+\!1)^{1\!/\!2+\epsilon}} {\Vert\nabla f(\mathbf{x}_t)\Vert}^2 (\omega) \le C(\omega) \ \textnormal{a.s.}
\end{equation}
By multiplying both the numerator and denominator of $t$-th term on the left side by $(t\!+\!1)^{1\!/\!2-\epsilon}$, the above inequality can be written as
\begin{equation}
\sum_{t=0}^\infty \frac{1}{t\!+\!1} [(t\!+\!1)^{1\!/\!2-\epsilon} {\Vert\nabla f(\mathbf{x}_t)\Vert}^2](\omega) \le C(\omega) \ \textnormal{a.s.}
\end{equation}
\par
Proof by contradiction leads to the desired result.
\end{proof}
\par
\vspace{2pt}
The stepsize is not required to be nonincreasing in our analysis.
But the finite $l^{1+\alpha}$ norm of the stepsize is necessary, which is easier to fulfill in practice.
It is worth noting that the assumption on boundness of the gradient $\nabla f$ and the oracle $\mathbf{g}$ is extremely common in articles analyzing adaptive methods, though it seems somewhat unsatisfying from a theoretical point of view.
The independence of $\gamma_t$ and $\xi_t$ is also of importance.
If $v_t = \beta'v_{t-1} + (1\!-\!\beta'){\Vert\mathbf{g}_t\Vert}^2$, the conditional expectation of update direction $\gamma_t\mathbf{g}_t$ will not be guaranteed to make an acute angle with the accurate gradient.
In this case, we need additional conditions, e.g., the limit on oscillation of effective stepsizes, which exceeds the scope of this article.

\section{Conclusion and Discussion}
We have provided mild conditions to ensure $L^2$ and almost sure convergence of stochastic gradient algorithms with momentum terms in the nonconvex setting.
The analysis is presented within a general framework while some common assumptions are weakened in this paper.
Particularly, $\nabla f$ is permitted to be $\alpha$-H{\"o}lder continuous.
Moreover, we go in the direction of showing the convergence of a modified version of AdaGrad and Adam.
\par
Similar extensions to original adaptive algorithms, however, are more complicated.
Our current analysis here does not necessarily hold in the case, mainly because the stepsize $\gamma_t$ is a function of past gradients $\mathbf{g}_0, \dots, \mathbf{g}_t$ and the expected update direction may deviate from the exact gradient.
Another limitation is the fact that the paper provides limited guidance on how to set the parameters such as stepsizes and the momentum factor in practice.
We leave these possible extensions as interesting topics for future research.

\newpage
\printbibliography

\end{document}